\documentclass[twoside,draft,reqno]{birkart}

\usepackage{amssymb}

\newcommand{\Vlad}{{Vl\u{a}du\c{t}}}
\newcommand{\DV}{{Drinfeld-$\!$\Vlad}}
\def\be{\begin{equation}}
\def\ee{\end{equation}}
\def\bea{\begin{eqnarray}}
\def\eea{\end{eqnarray}}
\def\Ga{{\bf G}_a}
\def\Ltau{{L\{\tau\}}}
\def\Lbartau{{\bar L\{\tau\}}}

\def\X{{\rm X}}
\def\Xo{\X_0}
\def\dX{{\dot\X}}
\def\dXo{{{\dot\X}_0}}
\def\GL{{\rm GL}}
\def\0{^{\phantom0}}
\def\9{_{\phantom9}}
\def\ra{\rightarrow}

\def\geqs{\geqslant}
\def\leqs{\leqslant}

\begin{document}

\title[Explicit towers of Drinfeld modular curves]
{Explicit towers of Drinfeld modular curves}
\author[N.D.~Elkies]{Noam D. Elkies}

\address{
Department of Mathematics,\br
Harvard University,\br
Cambridge, MA 02138 USA}

\email{elkies@math.harvard.edu}

\begin{abstract}
We give explicit equations for the simplest towers of Drinfeld modular
curves over any finite field, and observe that they coincide with
the asymptotically optimal towers of curves constructed by
Garcia and Stichtenoth.
\end{abstract}

\maketitle

\section{Introduction}  Fix a finite field~$k_1$ of size $q_1$.
It has been known for almost twenty years (see \cite{Ihara,TVZ,TV}) that
any curve $C/k_1$ of genus~$g=g(C)$ has at most $(q_1^{1/2}-1+o(1))g$
points rational over~$k_1$ as $g\ra\infty$
(the {\em \DV\ bound}\/~\cite{DV}),
and that if $q_1$ is a square then there are various families of
classical, Shimura, or Drinfeld modular curves~$C/k_1$ of genus
$g\ra\infty$ with $\#(C(k_1)) \geq (q_1^{1/2}-1)(g-1)$.
Thus such curves are ``asymptotically optimal'': they attain the
lim~sup of $\#(C(k_1))/g(C)$.  Moreover, asymptotically optimal curves
yield excellent linear error-correcting codes over~$k_1$~\cite{Goppa,TV}.

To actually construct and use these Goppa codes one needs explicit
equations for asymptotically optimal~$C$.  Now the definitions of
modular curves are in principle constructive, but it is usually not
feasible to actually exhibit a given modular curve.  For certain
Drinfeld modular curves, \cite[pp.453ff.]{TV} gives explicit,
albeit unpleasant, models as plane curves with two complicated
singularities.  Garcia and Stichtenoth gave nice formulas
for asymptotically optimal families of curves
in~\cite{GS:classical} for $q_1=4,9$ (see also \cite{GST}), and
in~\cite{GS:drinfeld} for all square prime powers~$q_1$.
In each case they construct a sequence of curves $C_1,C_2,C_3,\ldots$
forming what we'll call a {\em recursive tower.}  A ``tower''
$\{C_n\}$ is a sequence of curves in which each $C_n$ is given
as a low-degree cover of $C_{n-1}$.  We shall call a tower
``recursive'' if $C_2$ is given as a curve in $C_1\times C_1$,
and $C_n$ is the curve in $C_1^n$ consisting of $n$-tuples
$(P_1,\ldots,P_n)$ of points in~$C_1$ such that $(P_j,P_{j+1})\in C_2$
for each $j=1,2,\ldots,n-1$.  That is, $C_n$ is obtained by iterating
$n-1$ times the correspondence from~$C_1$ to itself given by~$C_2$.
In each of the cases discovered by Garcia and Stichtenoth,
analysis of this correspondence yields the genus of each $C_n$
and enough $k_1$-rational points on~$C_n$ to prove that
$\#(C_n(k_1))\geq (q_1^{1/2}-1) g(C_n)$.
The genus grows exponentially with~$n$, yet only $O(n)$ equations
of bounded degree are needed to exhibit $C_n$.

Now while equations for most modular curves cannot be reasonably
exhibited, modular curves whose conductors are products of small primes
form towers, and modular curves with a single repeated factor
even form recursive towers.
For instance, for any integers $N_0\geq1$ and $l>1$, the
classical modular curves $X_0(N_0 l^n)$ ($n=1,2,3,\ldots$) form
a recursive tower of curves related by maps of degree~$l$.
There are analogous towers of Shimura and Drinfeld modular curves.
We observed in~\cite{NDE:Allerton} that this can be used to exhibit
asymptotically optimal towers, and carried out the computation of
explicit equations for eight such towers, six of classical modular
curves and two of Shimura curves.  Moreover, we showed that the
towers of classical modular curves $\Xo(2^n)$ and $\Xo(3^n)$
in characteristics~$3$ and~$2$ respectively are the same as
the towers constructed in~\cite{GS:classical}.

In the same paper we noted that similar methods can be used to
exhibit towers of Drinfeld modular curves, and announced that
one such tower recovers the curves constructed in~\cite{GS:drinfeld}
for arbitrary square~$q_1$.  In the present paper we specify this tower
and perform the computations that determine its equations and thus
identify it with the Garcia-Stichtenoth tower.  We do the same for
a closely related tower obtained by Garcia and Stichtenoth
in~\cite{GS:more_drinfeld}.

In the next section we recall basic definitions of Drinfeld modules,
isogenies, and super\-singularity.  The following section describes
certain Drinfeld modular curves.  In both sections we simplify the
exposition by describing only the modules and curves that arise in
the interpretation of the Garcia-Stichtenoth towers; for a thorough
treatment of the general case we refer to \cite{Gekeler:LNM}.  In the
final section we give explicit equations for several of the simplest
towers of Drinfeld modular curves and observe that two of them
coincide with asymptotically optimal towers obtained by
Garcia and Stichtenoth.

{\bf Acknowledgements.}  I thank Bjorn Poonen for introducing me
to Drinfeld modules, and Henning Stichtenoth for information
on his and Garcia's asymptotically optimal towers.  Thanks also
to E.-U.~Gekeler, D.~Goss, B.H.~Gross, B.~Poonen, and M.~Zieve
for much enlightening communication concerning Drinfeld modules
and modular curves.

This work was made possible in part by funding from the
Packard Foundation.

\section{Drinfeld modules and isogenies}
Fix a finite field~$k$ of size~$q$.
(In our application, $k_1$ will be the quadratic extension of~$k$,
so $q\0_1=q^2$; the use of ``$k$'' instead of ``${\bf F}_q$'' is
our only divergence from the notations of \cite{Gekeler:LNM}).
For any field $L\supseteq k$, we denote by $\Ltau$
the non-commutative $L$-algebra generated by $\tau$ and satisfying
the relation $\tau a = a^q \tau$ for all~$a\in L$.  Equivalently,
$\Ltau$ is the ring of endomorphisms of $\Ga$ defined over~$L$ and
linear over~$k$.  Explicitly, the polynomial $\sum_{i=0}^n l_i \tau^i$
acts on~$\Ga$ as the endomorphism taking any~$X$\/ to the $q$-linearized
polynomial $\sum_{i=0}^n l_i X\9^{q^i}$.

We only consider Drinfeld modules of rank~$2$, and usually only ones
associated to a function field of genus~$0$ with a place at infinity
of degree~$1$.  Let then $K=k(T)$ and $A=k[T]$.  (See
\cite{Gekeler:LNM} for the general case, in which $A$ is the ring of
functions on a curve over~$k$ with poles at most at a fixed place
$\infty$ of the curve.)  In general, a Drinfeld module is defined as a
$k$-algebra homomorphism $\phi: a \mapsto \phi_a$ from~$A$ to~$\Ltau$
satisfying certain technical conditions.  In our case, $A=k[T]$,
so specifying $\phi$ is equivalent to choosing $\phi\0_T$.
The rank of the resulting Drinfeld module is then simply the
degree of $\phi\0_T$ as a polynomial in~$\tau$.  Thus for us
\begin{equation}
\phi\0_T = l_0 + l_1 \tau + l_2 \tau^2 = l_0 + g\tau + \Delta\tau^2,
\label{phi(T)}
\end{equation}
with nonzero {\em discriminant}~$\Delta=\Delta(\phi)$.
In general the map~$\gamma: A \ra L$ taking any $a\in A$ to the
``constant term'' ($\tau^0$ coefficient) of $\phi_a$ is a ring
homomorphism; in our case $\gamma$ is determined by $\gamma(T)=l_0$.

If $\phi,\psi$ are two Drinfeld modules, an
{\em isogeny} from $\phi$ to~$\psi$ is a $u\in\Lbartau$ such that
\begin{equation}
u \circ \phi_a = \psi_a \circ u
\label{isogeny}
\end{equation}
for all $a\in A$.  For our $A$, (\ref{isogeny}) holds for all $a\in A$
if and only if it holds for~$a=T$.  The {\em kernel} of the isogeny
is\footnote{
  When $u$ is not separable, i.e.\ has $\tau^0$ coefficient zero,
  it is for many purposes better to consider $\ker(u)$ not as a
  subgroup of $\bar L$ but as a group subscheme of~$\Ga$.  We shall
  not need this refinement here.  Note that the condition
  $\gamma(a)\neq0$ occurring later in this paragraph is
  equivalent to the separability of $\phi_a$.
  }
\begin{equation}
\ker(u) := \{ x \in \bar L: \phi_u(x) = 0\}.
\label{ker(u)}
\end{equation}
This is a $k$-vector subspace of~$\bar L$, which is of finite
dimension unless $u=0$.  By (\ref{isogeny}), $\phi_a(x)\in\ker(u)$
for all $x\in\ker(u)$, so $\ker(u)$ in fact has the structure of an
$A$-module.  Conversely, for every finite $G\subset\bar L$ which
is an $A$-submodule of~$\bar L$ for the $\phi$-action of~$A$
on~$\bar L$, one may define $u\in\Lbartau$ of degree $\dim_k G$\/ by
\begin{equation}
u(X) = \prod_{x\in G} (X-x),
\label{udef}
\end{equation}
and then $u$ is an isogeny with kernel~$G$\/ from~$\phi$ to some
Drinfeld module~$\psi_a$.  In particular, if $u=\phi_{a_1}$ then
(\ref{isogeny}) holds with $\psi=\phi$ for any $a_1\in A$; thus
$\phi_{a_1}$ is an isogeny from~$\phi$ to itself,
called {\em multiplication by~$a_1$}.
If $\gamma(a_1)\neq 0$, the kernel of this isogeny is isomorphic with
$(A/a_1 A)^2$ as an $A$-module \cite[Prop.~I.1.6]{Gekeler:LNM}.
Elements of $\ker(\phi_{a_1})$ are called {\em $a_1$-division points}
or {\em $a_1$-torsion points} of the Drinfeld module~$\phi$.
In particular, the $T$-torsion points are the roots in~$x\in\bar L$ of
\begin{equation}
\phi\0_T(X) = \gamma(T) X + g X^q + \Delta X^{q^2}.
\label{T-torsion}
\end{equation}
If $\gamma$ is not injective then $\ker\gamma = A a_0$ for some
irreducible $a_0\in A$.  Then the \hbox{$a_0$-torsion} points of~$\phi$
constitute a vector space of dimension $1$ or~$0$ over the field
$A/a_0 A$.  The Drinfeld module~$\phi$ is then said to be
{\em supersingular} if $\ker(\phi_{a_0})=\{0\}$, {\em ordinary}
otherwise.  We shall use the case $\deg(a_0)=1$, when
$\phi_{a_0}=g\tau+\Delta\tau^2$ and thus $\phi$ is supersingular
if and only if $g=0$.  Note that $a_0$ is of degree~$1$ if and
only if $\gamma(T)\in k$, and that Drinfeld modules over~$\bar k$
with $\gamma(T)\in k$ may arise as the ``reduction mod~$(T-\gamma(T))$''
of Drinfeld modules over~$\bar k(T)$ with $\gamma(T)=T$.

\section{Drinfeld modular curves}

An {\em isomorphism} between Drinfeld modules is an invertible
isogeny, i.e.\ some $u\in\bar L^*$ satisfying~(\ref{isogeny}).
This isomorphism multiplies each coefficient $l_i$ in~(\ref{phi(T)})
by $u^{1-q^i}$.  We define the {\em $J$-invariant} of a Drinfeld
module~$\phi$ given by~(\ref{phi(T)}) as follows:\footnote{
  Usually a lower-case $j$ is used for this.  We use a capital~$J$\/
  to forestall confusion with the integer variable $j$ appearing in
  the next section.  For elliptic curves one sometimes sees $J=j/12^3$;
  in the Drinfeld modular setting, no factor analogous to $12^3$
  is needed, so we might plausibly claim that the invariant
  of a Drinfeld module corresponds to~$J$\/ as well as~$j$ in the
  classical theory of modular invariants of elliptic curves\ldots
  }
\begin{equation}
J(\phi) = \frac{g^{q+1}}{\Delta}.
\label{J(phi)}
\end{equation}
Two Drinfeld modules with the same~$\gamma$ are isomorphic
(over~$\bar L$) if and only if their $J$-invariants are equal.
Thus, in analogy with the case of the classical modular curves
parametrizing elliptic curves, we refer to the $J$-line as the
{\em Drinfeld modular curve} $\X(1)$ for Drinfeld modules with
a given $\gamma$.  Likewise, for $N\in A$\/ such that $\gamma(N)\neq0$,
we have Drinfeld modular curves $\X_0(N)$, $\X_1(N)$, $\X(N)$
parametrizing Drinfeld modules with a given $\gamma$ and a choice
of torsion subgroup $G \cong A/NA$, or such a subgroup~$G$\/
together with a generator of~$G$\/ as an $A$-module, or an
identification of the group of $N$\/-torsion points with $(A/NA)^2$.
These are finite separable covers of~$\X(1)$, all of which except
$\X(N)$ ($N\notin k^*$) are geometrically irreducible; $\X(N)$ is a
normal cover of~$\X(1)$ with Galois group $\GL_2(A/NA)$, and $\X_1(N)$
is an abelian normal cover of $\X_0(N)$ with Galois group $(A/NA)^*$.
Note that, unlike $\X(1)$, these curves $\X_0(N)$, $\X_1(N)$, $\X(N)$
generally depend on the choice of $\gamma$.  If $\gamma(T)\in k$,
we may regard the curves $\X(1),\X_0(N),\X_1(N),\X(N)$
as the ``reduction mod~$(T-\gamma(T))$'' of the corresponding modular
curves for $\gamma(T)=T$.  More generally, reducing the $\gamma(T)=T$\/
curves modulo any irreducible $a_0\in A$ yields the curves parametrizing
Drinfeld modules for which $\gamma(T)$ is a root of~$a_0$.

If $\gamma$ is not injective, we say that a point on a
Drinfeld modular curve is {\em ordinary} or {\em supersingular}
according as the Drinfeld modules it parametrizes are ordinary
or supersingular.  It is known that in this case the supersingular
points constitute a nonempty finite set.  For instance, we have
seen in effect that if $\gamma(T)\in k$ then $\X(1)$ has the unique
supersingular point $J=0$.  The supersingular points on $\X_0(N)$,
$\X_1(N)$ and $\X(N)$ are then the preimages of $J=0$ under the natural
maps from those Drinfeld modular curves to~$\X(1)$.  It is known that
each supersingular point on~$\X_0(N)$ is defined over the quadratic
extension~$k_1$ of~$k$, and that the same is true for certain
twists of $\X_1(N)$ and of each component of~$\X(N)$.

We shall relate the Garcia-Stichtenoth curves to certain Drinfeld
modular curves with $N=T^n$ and $\gamma(T)=1$.  We shall find that
in some cases the modules parametrized by these curves must satisfy
the additional condition $\Delta=-1$, i.e.\
\begin{equation}
\phi\0_T = 1 + g\tau - \tau^2.
\label{normalized}
\end{equation}
We call such Drinfeld modules {\em normalized}.  A Drinfeld
module~$\phi$ with $\gamma(T)=1$ is isomorphic to a normalized one
if and only $-\Delta(\phi)$ is a $(q^2-1)$-st power.  This condition
is invariant under isogeny.  One thus expects that there would be an
equivalent condition in terms of the torsion structure of~$\phi$.
Such a condition cannot be given in completely elementary terms,
because $\GL_2(A)$ does not have a large enough cyclic quotient.
But\footnote{
  Thanks to Bjorn Poonen for pointing this out, and to David Goss for
  the reference to Hamahata.
  }
one can give an equivalent condition in terms of $\wedge^2\phi$.
For a general Drinfeld module~$\phi$ given by (\ref{phi(T)}), 
``$\wedge^2\phi$'' is the rank-$1$ module $T\mapsto l_0-l_2\tau$,
whose Tate modules are the discriminants of those of $\phi$
\cite[Thm.4.1]{Hamahata}.  Thus $\phi$ is normalized 
if and only if $\wedge^2\phi$ takes $T$\/ to $1-\tau$.

In terms of the coordinate~$J$\/ on $X(1)$, the condition that
$\phi$ be equivalent to a normalized module is that
$-J(\phi)$ be a $(q+1)$-st power.  Thus normalized modules,
or normalized modules with suitable level-$N$\/ structure,
are parametrized by curves we shall call $\dX(1)$, $\dXo(N)$,
$\dX_1(N)$, $\dX(N)$, whose function fields are obtained from those of
$\X(1),\Xo(N),\X_1(N),\X(N)$ by adjoining a $(q+1)$-st root of $(-J)$.

Now by analogy with the case of classical and Shimura curves, one
might expect that the curves $\Xo(N)$ are asymptotically optimal
over $k_1$, and that the same is true of the twists mentioned earlier
of $\X_1(N)$ and of the components of $\X(N)$,
with supersingular points already providing
$(q-1+o(1))g$ rational points in each case.  One might
even hope that the same is true with $\Xo$, $\X_1$, $\X$ replaced
by $\dXo$, $\dX_1$, $\dX$.  It would be enough to prove this for
$\dX(N)$, since all the other curves listed are quotients of $\dX(N)$
by subgroups of $\GL_2(A/NA)$.  This is stated explicitly in
the literature only for components of $\X(N)$ in the case that
$N$\/ is an irreducible polynomial of odd degree (see for instance
\cite[pp.449ff.]{TV}).  We expect that the same is true for
$\dX(N)$ and arbitrary~$N$, and indeed even for modular curves
for Drinfeld modules over rings other than $k[T]$.  In each case the
supersingular points are readily enumerated, and the main technical
challenge is computing the genus of the curve, since the covering maps
to $\X(1)$ are highly and wildly ramified above the ``cusp'' $J=\infty$.
[This was also the most difficult part of Garcia and Stichtenoth's
direct construction in~\cite{GS:drinfeld}.]

Even though explicit statements of asymptotic optimality have not
been made, the genera of $\Xo(N)$, $\X_1(N)$, and the components
of~$\X(N)$ are known.  They were computed in \cite[Thm.4.4]{Goss}
for~$\X(N)$, and also in Gekeler's
thesis~\cite[Satz~3.4.8]{Gekeler:thesis},
which also deals with the case of $\Xo(N)$ (Satz~3.4.18).
The curves $\X_1(N)$ were treated in~\cite{GekNon}.
These results, combined with the enumeration of supersingular points,
should yield the asymptotic optimality of all these curves over~$k_1$.
Alternatively, M.~Zieve suggests, and Gekeler confirms by e-mail,
that one may be able to entirely avoid the genus computation and the
enumeration of supersingular points by adapting an earlier proof
by Ihara \cite[pp.292--3]{Ihara:Corvallis}
that a classical modular curve of genus~$g$ has at least $(p-1)(g-1)$
points rational over the field of $p^2$ elements.  That argument uses
the reduction mod~$p$ of certain Hecke correspondences on the curve,
such as $\X_0(Np)$ considered as a correspondence on $\X_0(N)$ via
its map to $\X_0(N)\times \X_0(N)$ (when $\gcd(p,N)=1$); similar
correspondences are available in the Drinfeld modular setting.
We are not aware of a published analysis along the same lines of
the curves $\dX_0(N)$ etc.; but as Gekeler points out it should be
straightforward to obtain their genera from those of $\X_0(N)$ etc.,
of which they are cyclic covers of degree prime to~$q$.
At any rate the results of the next section, together with the
genus calculations in~\cite{GS:drinfeld,GS:more_drinfeld},
show at least that the curves $\dX_0(T^n)$ are asymptotically optimal.

\section{Some Drinfeld modular curves of conductor $T^n$}

If $x_1$ is a nonzero torsion point of~$\phi$ then we can solve
for~$g$ by setting (\ref{T-torsion}) equal to zero, obtaining
\begin{equation}
0 = T(x_1) = x_1 + gx_1^q - x_1^{q^2}, \quad
{\rm i.e.} \quad
g = x_1^{-q} (x_1^{q^2}-x_1).
\label{x1}
\end{equation}
Thus $x_1$ may be regarded as a coordinate for the rational curve
$\dX_1(T)$ parametrizing normalized Drinfeld modules $\phi$ with a
$T$\/-torsion point. For any nonzero~$x\in\bar L$ we let $t_x$
be the corresponding linearized polynomial
\begin{equation}
t_x(X) := X + \frac{x^{q^2}-x}{x^q} X\9^q - X\9^{q^2}.
\label{t_x1}
\end{equation}
The supersingular $x_1$ are those for which the $X\9^q$ coefficient
of $t_{x_1}$ vanishes, i.e.\ the units of the quadratic extension~$k_1$
of~$k$.  (The points $x_1=0,\infty$ are the cusps of $\dX_1(T)$.)
We next construct the curves we shall call $\dX'_0(T^n)$, parametrizing
$\phi$ with such a torsion point $x_1$ as well as a torsion group
$G_n\cong A_0[T]/T^n$ containing $x_1$, and find the supersingular
points on these curves.  Note that $\dX'_0(T)=\dX_1(T)$, but for
$n\geqs2$ the curve $\dX'_0(T^n)$ is only a quotient of $\dX_1(T^n)$,
because, as in the case of elliptic curves, we demand not that each
point of $G_n$ be rational, only that $G_n$ be permuted by Galois.

For $j\leqs n$ let $G_j$ be the group $T^{n-j} G_n$ of $T^j$-torsion
points in $G_n$.  Of course $G_1 = k x_1$.  For any $x\neq0$ let $P_x$
be the linearized polynomial
\begin{equation}
P_x(X) = x^{q-1} X - X\9^q
\label{Px}
\end{equation}
vanishing on $kx$.  Since $t_x$ vanishes on $kx$ we
expect it to factor through $P_x$ (see for instance
\cite[Prop.3]{NDE:qalg}), and indeed we find
\begin{equation}
t_x(X) = Q_x(P_x(X))
\label{t=QP}
\end{equation}
where $Q_x$ is the linearized polynomial
\begin{equation}
Q_x(X) = x^{1-q} X + X\9^q.
\label{Qx}
\end{equation}
Let $t'_x(X)$ be the reverse composition $P_x \circ Q_x$, given by
\begin{equation}
t'_x(X) := P_x(Q_x(X)) =
X + (x^{q-1}-x^{q-q^2}) X\9^q - X\9^{q^2}.
\label{PQ}
\end{equation}
Note that $t'_x$ again gives a normalized Drinfeld module of rank~2,
namely the module $T$\/-isogenous\footnote{
  i.e.\ with an isogeny (here $P_x$ or $Q_x$)
  whose kernel is isomorphic with $A/TA$ as an $A$-module.
  }
to $\phi$ obtained as the quotient of $\phi$ by $G_1$.
Now suppose $y_2$ is a generator of $G_2$ such that
$t_{x_1}(y_2) = x_1$.  There are $q$ such $y_2$, all
differing by multiples of $x_1$, so $x_2 := P_{x_1}(y_2)$
does not depend on the choice of $t_2$; conversely $x_2$
determines $G_2$.  Since $t_{x_1} = Q_{x_1} \circ P_{x_1}$,
the condition $t_{x_1}(y_2)=x_1$ is equivalent to
$Q_{x_1} (x_2) = x_1$.  Recalling the definition (\ref{Qx})
and multiplying by $x_1^q$ we find
\begin{equation}
x_2 = x_1^{-1} z_2 \quad {\rm where} \quad z_2^q + z_2 = x_1^{q+1}.
\label{Fermat!}
\end{equation}
Thus the curve $\dX'_0(T^2)$ is just $z_2^q + z_2 = x_1^{q+1}$,
which is known to be $k_1$-isomorphic with the Fermat curve
of degree $q+1$ (a.k.a.\ the ``Hermitian curve'' over~$k_1$).
Note that $G_2$ is the $k$-vector space of zeros of the linearized
polynomial $P_{x_2} \circ P_{x_1}$.  Moreover,
\begin{equation}
0 = P_{x_1}(Q_{x_1}(x_2)) = t'_{x_1}(x_2).
\label{t'=t2}
\end{equation}
That is, $x_2$ is a $T$\/-torsion point on our $T$\/-isogenous module.
It follows that \hbox{$t'_{x_1}=t_{x_2}$}.
By induction we can now determine the tower of curves $\dX'_0(T^n)$
explicitly: the function field of $\dX'_0(T^n)$
is generated by $x_1,x_2,\ldots,x_{n-1}$ with relations
$Q_{x_{j-1}} (x_j) = x_{j-1}$, or equivalently
\begin{equation}
x_j = x_{j-1}^{-1} z_j \quad {\rm where} \quad
z_j^q + z_j = x_{j-1}^{q+1}
\label{GSeq}
\end{equation}
($1<j<n$); the point $(x_1,x_2,\ldots,x_{n-1})$ parametrizes the
Drinfeld module with $\phi\0_T = t_{x_1}$, with $T^n$-torsion subgroup
generated by any of the $q^{n-1}$ solutions of
\begin{equation}
(P_{x_{n-1}} \circ P_{x_{n-2}} \circ \cdots \circ P_{x_2} \circ P_{x_1})
(y_n) = x_n.
\label{tn}
\end{equation}
This works because by applying
\begin{equation}
Q_{x_j} \circ P_{x_j} = T_{x_j} = P_{x_{j-1}} \circ Q_{x_{j-1}}
\label{dosido}
\end{equation}
$(n-2)$ times we find
\bea
x_{n-1} = Q_{x_{n-1}}(x_n) &=& (Q_{x_{n-1}} \circ P_{x_{n-1}}
\circ P_{x_{n-2}} \circ \cdots \circ P_{x_1}) (y_n)
\nonumber \\
&=&
(P_{x_{n-2}} \circ Q_{x_{n-2}} \circ P_{x_{n-2}}
\circ \cdots \circ P_{x_1}) (y_n)
\nonumber \\
= \ \cdots &=&
(P_{x_{n-2}} \circ \cdots \circ P_{x_2}
\circ Q_{x_1} \circ P_{x_1}) (y_n)
\label{gauntlet} \\
&=&
(P_{x_{n-2}} \circ \cdots \circ P_{x_2} \circ P_{x_1} \circ T_{x_1})
(y_n)
\nonumber \\
&=&
(P_{x_{n-2}} \circ \cdots \circ P_{x_2} \circ P_{x_1})
(T_{x_1} (y_n)),
\nonumber
\eea
i.e.\ $T_{x_1}(y_n)$ satisfies the equation for $y_{n-1}$.
As with $G_2$ we see that $G_n$ is the $k$-vector space of
zeros of the linearized polynomial
\begin{equation}
P_{x_n} \circ P_{x_{n-1}} \circ \cdots \circ P_{x_2} \circ P_{x_1}.
\label{Gn}
\end{equation}

In~\cite{GS:drinfeld} Garcia and Stichtenoth obtained many
$k_1$-rational points on these curves as follows.  For each
of the $q^2-1$ points on $\dX'_0(T)$ with $x_1\in k_1^*$ we have
$x_1^{q+1}\in k_1^*$, so the $q$ solutions $z_2$ of (\ref{Fermat!})
are just the $q$ elements of~$k_1$ whose trace to~$k$ is $x_1^{q+1}$.
Clearly none of these $z_2$ vanish, so $x_2 = x_1 z_2$ is again in
$k_1^*$.  Inductively we see that $x_1$ lies under $q^{n-1}$ rational
points of $\dX'_0(T^n)$ defined over~$k_1$.  Garcia and Stichtenoth
use the resulting $(q^2-1)q^{n-1}$ rational points to confirm that
the $\dX'_0(T^n)$ form an asymptotically optimal tower over~$k_1$.
From the Drinfeld modular viewpoint we recognize $x_1\in k_1^*$
as the condition that a point on~$\dX'_0(T^n)$ be supersingular.
Thus the $(q^2-1)q^{n-1}$ points of $\dX'_0(T^n)$ lying above
$k_1^*$ are precisely the supersingular points on $\dX'_0(T^n)$.

The group $k_1^*$ acts on $\dX'_0(T^n)$ by
\begin{equation}
c(x_1,\ldots,x_n) =
(c x_1, \bar c x_2, c x_3, \bar c x_4, \ldots, c^{q^{n-1}} x_n)
\quad (c\in k_1^*, \bar c := c^q).
\label{action}
\end{equation}
If $c\in k^*$ then the automorphism (\ref{action}) preserves~$\phi$
(see (\ref{x1})) and each~$G_j$, but changes the generator $x_1$
of~$G_1$ to~$cx_1$.  Thus we recover $\dX_0(T^n)$ as the quotient of
$\dX'_0(T^n)$ by the action of~$k^*$.  The tower of curves $\dX_0(T^n)$
was obtained in this way (again without mention or use of Drinfeld
modules) in~\cite{GS:more_drinfeld}.  For arbitrary $c\in k_1^*$, the
automorphism (\ref{action}) multiplies $g$ by the $(q+1)$-st root
of unity $c/\bar c$.  This takes $\phi$ to a Drinfeld module~$\phi_c$
with the same $J$-invariant but a different choice of $(-J)^{1/(q+1)}$.
The isomorphism~$c$ from~$\phi$ to~$\phi_c$ respects our choice of
subgroups~$G_j$.  Thus the quotient of $\dX_0(T^n)$ by $k_1^*/k\9^*$,
or equivalently of $\dX'_0(T^n)$ by $k)1^*$, is the Drinfeld modular
curve $\X_0(T^n)$.

We next obtain an explicit description of $\{\X_0(T^n)\}$ as
a recursive tower.  As with several of the examples
in~\cite{NDE:Allerton}, it will be convenient to start
the tower at $n=2$.  [To start at $n=1$ we would have to use
$x_j^{q^2-1}$ as the $j$-th coordinate, and then all the supersingular
points would be on the normalization of the highly singular point
$(x_1^{q^2-1},\ldots,x_n^{q^2-1})=(1,\ldots,1)$ on the resulting curve.]
The curve $\X_0(T)$ is the quotient of the $x_1$-line $\dX'_0(T)$
by $k_1^*$, so has genus zero and coordinate $x_1^{q^2-1}$.  To 
obtain $\X_0(T^2)$, raise both sides of the equation (\ref{Fermat!})
for $\dX'_0(T^2)$ to the power $q-1$ to obtain
\begin{equation}
Z_2(1+Z_2)^{q-1} = x_1^{q^2-1},
\quad{\rm where}\ Z_2 := z_2^{q-1} = (x_1 x_2)^{q-1}.
\label{T^2}
\end{equation}
Now $Z_2$ is invariant under $k_1^*$, and (\ref{T^2}) shows
that the $Z_2$-line is a degree-$q$ cover of the $x_1^{q^2-1}$-line
$\X_0(T)$; since the cover $\X_0(T^2)/\X_0(T)$ is also of degree~$q$,
we conclude that $Z_2$ generates the function field of $\X_0(T^2)$.

Thus for $n\geq 2$ the function field of $\X_0(T^n)$ is generated by
\begin{equation}
Z_j := z_j^{q-1} \quad (2\leq j \leq n).
\label{Zj}
\end{equation}
For each $j=2,\ldots,n-1$, we have
$Z_{j+1}(1+Z_{j+1})^{q-1} = x_j^{q^2-1}$
as in (\ref{T^2}), which in turn equals
\begin{equation}
 (z_j/x_{j-1})^{q^2-1}
= Z_j^{q+1}/x_{j-1}^{q^2-1}
= Z_j^q/(1+Z_j)^{q-1}.
\label{Zj+1}
\end{equation}
Thus
\begin{equation}
Z_{j+1} (1+Z_{j+1})^{q-1} = Z_j^q/(1+Z_j)^{q-1}.
\label{Zrec}
\end{equation}
This gives $Z_{j+1}$ as an algebraic function of degree~$q$
in~$Z_j$ (and vice versa).  Thus the relations (\ref{Zrec})
for $j=2,\ldots,n-1$ determine the function field of $\X_0(T^n)$.
{}From our description of the supersingular points on $\dX'_0(T^n)$
we see that the $q^{n-1}$ supersingular points on $\X_0(T^n)$
are the points for which each $Z_j$ is in
\bea
&   &\{ Z \in k_1: Z^{q+1}=1, Z\neq -1 \} \nonumber\\
& = &\{ Z \in k_1: Z(1+Z)^{q-1} = 1 \} \label{Zss} \\
& = &\{ Z \in k_1: Z^q = (1+Z)^{q-1} \}. \nonumber \\
\eea

\end{document}